\documentclass[graybox]{svmult}

\usepackage{makeidx}         
\usepackage{graphicx}        
\usepackage{multicol}        
\usepackage[bottom]{footmisc}
\usepackage{newtxtext}  
\usepackage{caption}
\usepackage{xcolor}
\usepackage{bm}
\usepackage{orcidlink}  %
\usepackage{algorithm2e}
\usepackage[varvw]{newtxmath}       
\newcommand{\setR}[0]{\mathbb{R}}
\newcommand{\wJ}[0]{\widetilde{J}}
\newcommand{\xb}[0]{\bm{x}}
\newcommand{\yb}[0]{\bm{y}}
\newcommand{\zb}[0]{\bm{z}}

\usepackage{cancel}
\makeindex             
\begin{document}

\title*{Variable reduction as a nonlinear preconditioning approach for optimization problems}
\titlerunning{Nonlinear preconditioning of optimization problems} 
\author{Gabriele Ciaramella \inst{1}\orcidlink{0000-0002-5877-4426} and Tommaso Vanzan \inst{2}\orcidlink{0000-0001-7554-4692}}
\institute{Gabriele Ciaramella \at MOX, Dipartimento di Matematica, Politecnico di Milano, Italy, \email{gabriele.ciaramella@polimi.it}, Member of the Indam GNCS,
\and Tommaso Vanzan \at Dipartimento di Scienze Matematiche, Politecnico di Torino, Italy, \email{tommaso.vanzan@polito.it}, Member of the Indam GNCS.}
%
%
\maketitle

\abstract*{When considering an unconstrained minimization problem, a standard approach is to solve the optimality system with a Newton method possibly preconditioned by, e.g., nonlinear elimination.
In this contribution, we argue that nonlinear elimination could be used to reduce the number of optimization variables by artificially constraining them to satisfy a subset of the optimality conditions. Consequently, a reduced objective function is derived which can now be minimized with \textit{any} optimization algorithm. By choosing suitable variables to eliminate, the conditioning of the reduced optimization problem is largely improved.
We here focus in particular on a right preconditioned gradient descent and show theoretical and numerical results supporting the validity of the presented approach.}

\abstract{When considering an unconstrained minimization problem, a standard approach is to solve the optimality system with a Newton method possibly preconditioned by, e.g., nonlinear elimination.
In this contribution, we argue that nonlinear elimination could be used to reduce the number of optimization variables by artificially constraining them to satisfy a subset of the optimality conditions. Consequently, a reduced objective function is derived which can now be minimized with \textit{any} optimization algorithm. By choosing suitable variables to eliminate, the conditioning of the reduced optimization problem is largely improved.
We here focus in particular on a right preconditioned gradient descent and show theoretical and numerical results supporting the validity of the presented approach.}

\section{Introduction}
In this manuscript, we consider the unconstrained optimization problem
\begin{equation}\label{eq:P}
\min_{\zb\in \setR^n} J(\zb),
\end{equation}
where $J:\setR^n\rightarrow \setR$ is a twice continuously differentiable cost function bounded from below and whose minimizers satisfy the nonlinear optimality conditions
\begin{equation}\label{eq:optimality_condition}
    \nabla J(\zb)=0.
\end{equation}
The direct solution of \eqref{eq:optimality_condition} using the Newton method can, in general, be hard, and may suffer from both poor initial guesses and a high computational cost. For these reasons, several algorithms (e.g., first-order methods, SQP, Quasi-Newton methods, trust-region methods) have been developed. From a Domain Decomposition (DD) community point of view, it is instead natural to study nonlinear preconditioning strategies to improve the convergence of Newton's method.
This is one of the most active research areas within the community as shown by the much participated minisymposium on the topic at the last conference at KAUST and the numerous contributions to the last conference proceedings, see, e.g., \cite{Brenner,Kohler,Kothari,Kothari2}. As nicely summarized in \cite[\S 1.2]{Liu2024}, nonlinear preconditioning strategies to solve a general nonlinear system $F(\zb)=0$ can be divided into left approaches (see, e.g., \cite{ASPIN,RASPEN,SRASPEN}), which replace the original system with an equivalent one $G(F(\zb))=0$ having the same roots but that is easier to solve, and right approaches (see, e.g., \cite{lanzkron1996,Cai2011}) which instead change the variables of the system by solving $F(G(\yb))=0$, with $\zb=G(\yb)$.

The goal of this contribution is to present an alternative approach to use a right-preconditioning approach, called nonlinear elimination, within the optimization field. Specifically, we use nonlinear elimination not to precondition the optimality system \eqref{eq:optimality_condition}, but to derive a reduced objective function which can then be minimized with any preferred optimization algorithm. The approach is promising both in terms of the theoretical analysis that can be developed and of computational efficiency.

To better frame our discussion in a DD setting, throughout this manuscript, we assume that the variable $\zb$ can be suitably split into $\zb=(\xb,\yb)$, with $\xb\in \setR^{n_x}$, $\yb\in \setR^{n_y}$ with $n_x+n_y=n$. We will not delve here into how this decomposition is derived: for certain optimization problems, $\zb$ is naturally decomposed into subsets of variables with clear different roles (as in risk-adverse PDE-constrained optimization which originally motivated this work  \cite[\S 5]{CMG}), see also  \cite[\S 1.5]{Liu2024} for an overview of spatial-, fields-, and physics-based decompositions in nonlinear preconditioning. We will however discuss that the efficacy of the methods analyzed strongly depends on a good choice for the variable partition.

\section{Nonlinear elimination and a right preconditioned gradient descent}
Given the decomposition $\zb=(\xb,\yb)$, \eqref{eq:optimality_condition} can be reformulated as
\begin{equation}\label{eq:optimality_conditions_expanded}
        \nabla_{\xb} J(\xb,\yb) = 0, \quad
        \nabla_{\yb} J(\xb,\yb) = 0.
\end{equation}
Even though we consider an unconstrained optimization problem, we could artificially\footnote{In contrast with cases where one can often eliminate some variables using the given constraints, problem \eqref{eq:P} is unconstrained, and the elimination arises from the optimality conditions.} constrain a set of variables in terms of the others. 
Assuming that for every $\xb$ the equation $\nabla_{\yb} J(\xb,\yb)=0$ admits a unique solution $\yb$ and that $\nabla_{\yb\yb} J(\xb,\yb)$ is invertible for every $(\xb,\yb)$, we can nonlinearly eliminate the $\yb$ variable by considering the differentiable implicit map $h:\setR^{n_x}\rightarrow \setR^{n_y}$ such that $\nabla_{\yb} J(\xb,h(\xb))=0$. We are then left to solve the reduced nonlinear equation $\widetilde{F}(\xb):=\nabla_{\xb} J(\xb,h(\xb))=0$, and to do so we may use the Newton iteration
\begin{equation}\label{eq:Newton_nonlinear_elimination}
\xb^{k+1}=\xb^k- (J\widetilde{F}(\xb^k))^{-1}\nabla_{\xb} J(\xb^k,h(\xb^k)),
\end{equation}
where a straight calculation using implicit differentiation shows that 
\[J\widetilde{F}(\xb)=\nabla_{\xb\xb} J(\xb,h(\xb)) -\nabla_{\yb\xb} J(\xb,h(\xb))\nabla_{\yb\yb}J(\xb,h(\xb))^{-1}\nabla_{\xb\yb}J(\xb,h(\xb)).\]
Iteration \eqref{eq:Newton_nonlinear_elimination} can be possibly globalized (or damped) by a line-search which would require the repeated (possibly inexact) evaluation of $h(\cdot)$ (and thus the solution of the equations $\nabla_{\yb} J(\xb,h(\xb))=0$).
What we have just described coincides exactly with the well-known nonlinear elimination method (see,e.g, \cite{lanzkron1996,Cai2011}) applied to \eqref{eq:optimality_conditions_expanded}.

However, since our original task is to solve an optimization problem, we remark that the variable elimination could be performed not only on the optimality system \eqref{eq:optimality_conditions_expanded}, but directly on the objective function. In other words, we consider the reduced cost function $\wJ:\setR^{n_x}\rightarrow \setR$ defined as $\wJ(\xb):=J(\xb,h(\xb))$. While a Newton's method applied to the optimality condition of $\wJ$ is equivalent to the nonlinear elimination approach recalled above for \eqref{eq:optimality_conditions_expanded} \footnote{Notice that $\nabla \wJ(\xb)=\nabla_{\xb} J(\xb,h(\xb))+\nabla_{\yb} J(\xb,h(\xb))h^\prime(\xb)=\nabla_{\xb} J(\xb,h(\xb))=\widetilde{F}(\xb)$, since by definition $h$ is such that $\nabla_{\yb} J(\xb,h(\xb))=0$.}, we are now actually free to choose our preferred algorithm to minimize $\wJ$. In particular, we here focus on a Gradient Descent (GD) method applied to $\wJ$, summarized by Alg. \ref{Alg:rightgradient}, and that from now on we call \textit{right} Preconditioned Gradient Descent (PGD) method, since we are performing a change of variables that guarantees that both the cost function and its gradient are always evaluated (even within the line-search) on points that satisfy a subset of the optimality conditions in \eqref{eq:optimality_conditions_expanded}.
This is contrast with the most popular approach to precondition GD, that is to multiply the gradient by a suitable invertible matrix which rescales the descent directions. This latter approach can be interpreted as a left preconditioning of the gradient descent method.

\begin{center}
\begin{minipage}{0.90\textwidth}%
\begin{algorithm}[H]
   \SetAlgoLined
   \KwData{Initial guess $\xb^0$, tolerance Tol.}
   \KwResult{Stationary point $\xb^\star$.}
   Set $k = 0$\;
   \While{$\|\nabla \wJ(\xb^k)\|>{\rm Tol}$}{
        Compute a step length $t$ such that $\wJ(\xb^k-t\nabla \wJ(\xb^k))<\wJ(\xb^k)$\;
        Set $\xb^{k+1}=\xb^k-t \nabla \wJ(\xb^k)$\;
   }
\caption{Right preconditioned gradient descent}\label{Alg:rightgradient}
\end{algorithm}
\end{minipage}
\end{center}

\begin{figure}
    \centering
    \includegraphics[width=0.32\linewidth]{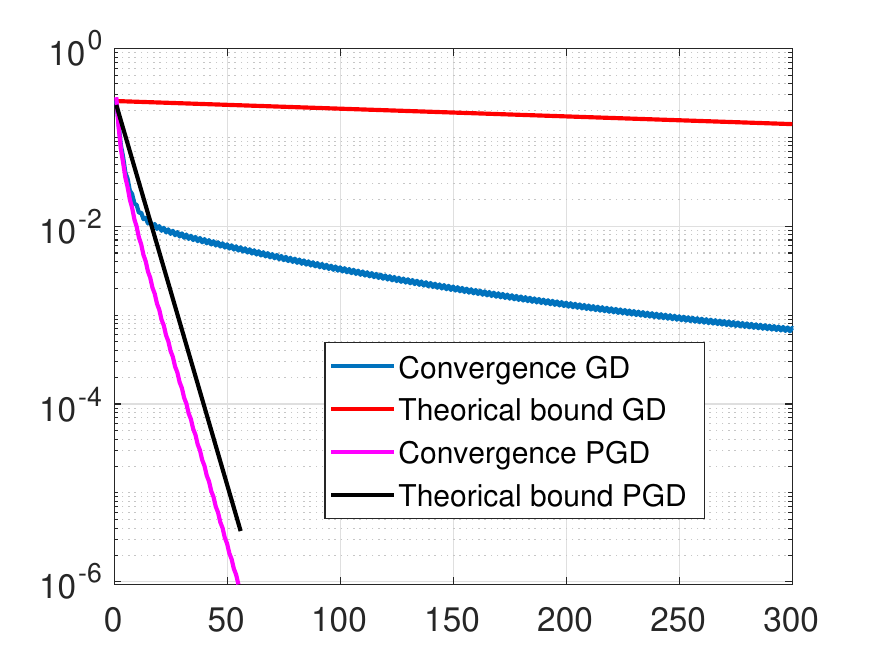}
    \includegraphics[width=0.32\linewidth]{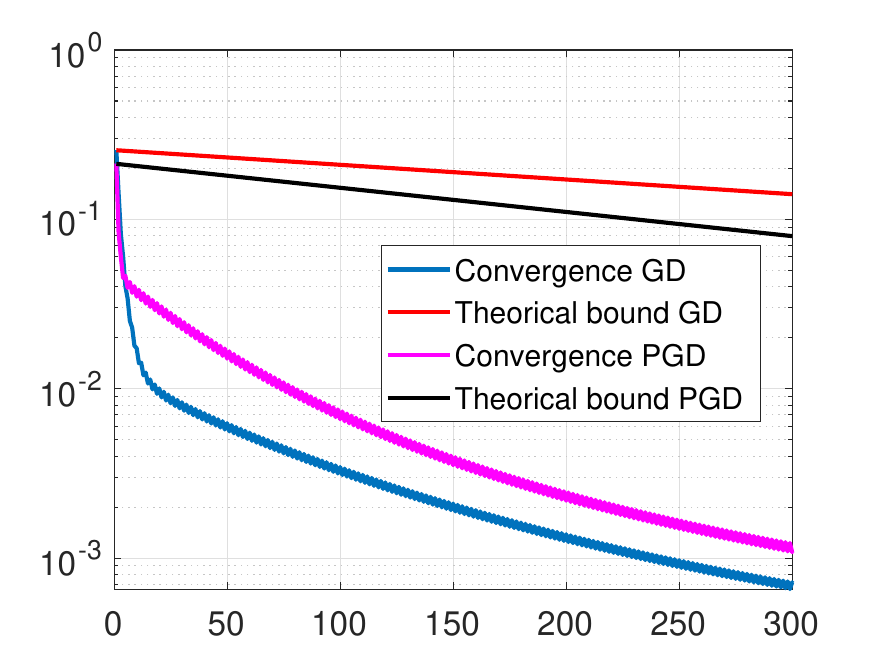}
    \includegraphics[width=0.32\linewidth]{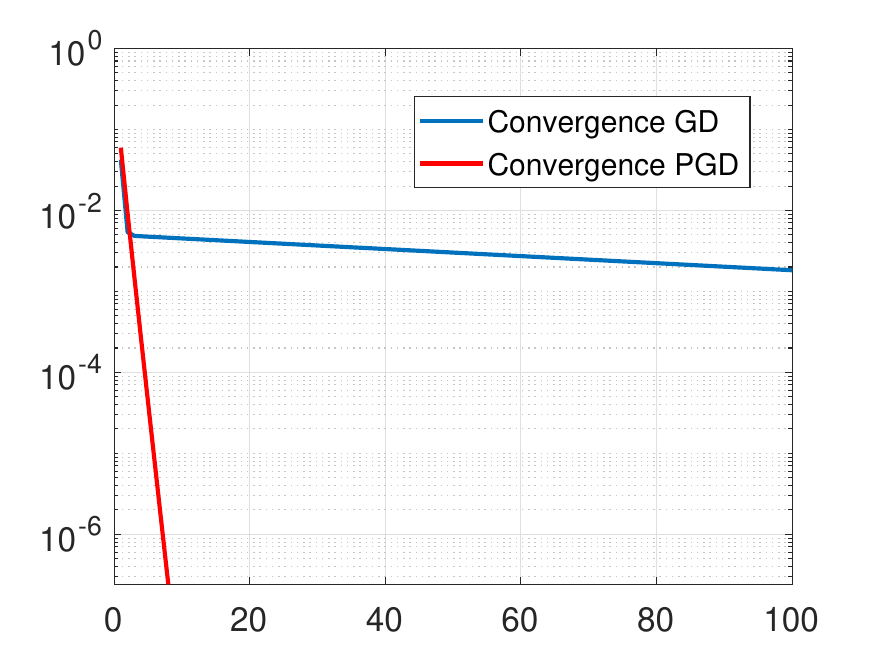}
    \caption{Comparison between the convergence of GD and PGD for a quadratic problem. $\kappa_2(A)= 1.001\cdot 10^3$, $\kappa_2(S)=10$ (left), $\kappa_2(S)= 6.08\cdot 10^2$ (center). Right panel convergence of GD and PGD for the strongly convex optimization problem \eqref{eq:strongly_convex} with $n=100$ and Armijo backtracking.}
    \label{fig:quadratic}
\end{figure}

In the rest of this section, we study the convergence properties of Alg. \ref{Alg:rightgradient} in a few relevant settings. As a first case study, we consider the quadratic cost function
\begin{equation}\label{eq:quadratic_problem}
J(\zb)= \frac{1}{2} \zb^\top A \zb -b^\top \zb +c=\frac{1}{2}\begin{pmatrix}\xb,\yb
\end{pmatrix}^\top \begin{pmatrix}
    A_{11} & A_{12}\\
    A_{21} & A_{22}
\end{pmatrix}\begin{pmatrix}\xb\\\yb
\end{pmatrix}-\begin{pmatrix}
    b_1 & b_2
\end{pmatrix}^\top \begin{pmatrix}\xb\\ \yb
\end{pmatrix} +c,
\end{equation}
$A$ being a s.p.d. matrix. It is well known (see, e.g., \cite[\S 1.3.2]{Bertsekas}) that for a quadratic problem there is an optimal choice for the step length (available in closed formula so that no line-search is needed) such that GD generates iterates satisfying 
\begin{equation}\label{eq:convergenge_GD}
\|\xb^k-\xb^\star\|_2\leq \sqrt{\kappa_2(A)}\left(\frac{\kappa_2(A)-1}{\kappa_2(A)+1}\right)^k\|\xb^0-\xb^\star\|_2.    
\end{equation}
The convergence rate then depends strongly on the conditioning of $A$. If we assume that $A_{11}$ is well conditioned while $A_{22}$ is ill conditioned, we may think about \textit{eliminating} the $\yb$ variables, deriving a reduced optimization problem only with respect to the $\xb$ variables, and hope to recover a fast convergence of GD applied to the reduced optimization problem.
To do so, we focus on the optimality condition of \eqref{eq:quadratic_problem}, corresponding to the linear system  $A\zb=\mathbf{b}$ from which, using static condensation, we may express $\yb=h(\xb)=-A_{22}^{-1}A_{21}\xb +A_{22}^{-1}\mathbf{b}_2$.
Inserting $\yb=h(\xb)$ into $J(\xb,\yb)$ we obtain the reduced cost function 
\begin{equation}\label{eq:reduced_cost_functional}
\wJ(\xb):=\frac{1}{2}\xb^\top S \xb+\widetilde{b}^\top \xb +\widetilde{c},
\end{equation}
where $S:=A_{11}-A_{12}A_{22}^{-1}A_{21}$ is the Schur complement and $\widetilde{b},\widetilde{c}$ are suitable vectors. Since the reduced cost function  $\wJ$ is still quadratic, the iterates of the right preconditioned GD (with the optimal step length) satisfy
\begin{equation}\label{eq:convergence_gradient_descent}
\|\xb^k-\xb^\star\|_2\leq \sqrt{\kappa_2(S)}\left(\frac{\kappa_2(S)-1}{\kappa_2(S)+1}\right)^k\|\xb^0-\xb^\star\|_2.
\end{equation}
Observing that
\begin{equation}\label{eq:argument_eigenvalues}
    \begin{aligned}
    \lambda_{\max}(S)&:=\max_{x\neq 0}\frac{\xb^\top S \xb}{\xb^\top \xb} \leq \max_{\xb\neq 0}\frac{\xb^\top A_{11}\xb}{\xb^\top \xb} \leq \max_{\zb=(\xb,\yb)\neq 0}\frac{\zb^\top A \zb}{\zb^\top \zb},\\
     \lambda_{\min}(A)&:=\min_{\zb=(\xb,\yb)\neq 0}\frac{\zb^\top A \zb}{\zb^\top \zb} \leq \min_{\xb:\zb=(\xb,-A_{22}^{-1}A_{21}\xb} \frac{\xb^\top S \xb}{\xb^\top \xb + \xb^\top A_{12}A_{22}^{-2}A_{21}\xb }\leq \lambda_{\min}(S) ,\\    
\end{aligned}
\end{equation}
we conclude that $\kappa_2(S)\leq \kappa_2(A)$, and we have prove the following proposition.
\begin{proposition}\label{prop:1}
For an unconstrained quadratic optimization problem, the right PGD method always has a better asymptotic convergence rate than the standard GD method.
\end{proposition}
While Proposition \ref{prop:1} guarantees that the right PGD converges asymptotically better than GD applied to the original problem \textit{regardless} of the variable decomposition, the latter still plays a key-role in determining the actual speed-up. To see this, we consider a matrix $A$ where $A_{11}$ and $A_{22}$ are similar through random orthogonal matrices to diagonal matrices with equispaced values between $10$ and $1$, and $10^3$ and $1$ respectively, and $A_{12}=A_{21}^\top$ is a small random perturbation such that $A$ is still s.p.d.. The matrix sizes are $n_x=40$ and $n_y=60$. On the left panel of Figure \ref{fig:quadratic}, we show the convergence of GD and of PGD where we fully eliminate the variable $\yb$. It is evident that ill conditioning due to the $\yb$ variable has completely vanished (this is confirmed by $\kappa_2(S)$ which turns equal to $\kappa_2(A_{11})$). Standard GD applied to $J(\xb,\yb)$ requires 3004 iterations and 0.12 seconds to reach a tolerance of $10^{-6}$ on the relative gradient norm. The PGD requires only 56 iterations and 0.03 seconds, thus it leads to a speed up of a factor $4$. Note that the $S$ is not assembled, but handled in a matrix-free way for a fair comparison, so that each evaluation of $h$ requires a linear solve.
On the center panel instead, we only eliminate the last $n_r=50$ variables of $\yb$. In this case, even though the theoretical bound for PGD is slightly better, there is no actual gain from the variable elimination: PGD needs now 2219 iterations, still less than GD, but the computational time rises up to 0.87 seconds . As a general rule, $h$ should satisfy the same properties of a linear preconditioner: it should be cheap to compute (even more important than in the linear case since $h$ is evaluated also during the line-search), and lead to a much better conditioned reduced problem. This is heuristically achieved by eliminating all variables associated to small eigenvalues (i.e. flat valleys) of $\nabla^2 J$ evaluated at the optimal value $\zb^\star$.

We now remark that Proposition \ref{prop:1} can be readily extended to
strongly convex cost functions that have Lipschitz continuous gradients (\cite[Chapter 1]{Nesterov2018}), namely functions that for certain $\mu,L\in \setR$ satisfy
\begin{equation}
    \begin{aligned}
    &J(\zb^\prime)\geq J(\zb)+\langle \nabla J(\zb),\zb^\prime-\zb\rangle +\frac{\mu}{2} \|\zb-\zb^\prime\|^2_2,\quad \forall \zb,\zb^\prime\in \setR^n,\\
    &\|\nabla J(\zb)-\nabla J(\zb^\prime)\|_2 \leq L \|\zb-\zb^\prime\|^2_2,\quad \forall \zb,\zb^\prime\in \setR^n.
    \end{aligned}
\end{equation}
For such functions, it is possible to prove (\cite[Theorem 2.1.15]{Nesterov2018}) that \eqref{eq:convergenge_GD} holds with $\kappa_2(A)$ replaced by $\frac{L}{\mu}$.
Since then 1) for twice continuously differentiable functions $\mu$ and $L$ can be taken equal to the minimum and maximum eigenvalues of $\nabla^2 J(\zb)$ over all $\zb\in \setR^n$, 2) $\nabla_{\xb\xb}^2 \wJ(\xb) $ corresponds to the Schur complement of $\nabla^2 J(\zb)$ for all  $\zb=(\xb,h(\xb))$ and thus the strongly convex parameters $\widetilde{\mu},\widetilde{L}$ of $\wJ$ satisfy (repeating the arguments of \eqref{eq:argument_eigenvalues} for all $\xb$) $\mu\leq \widetilde{\mu}$ and $\widetilde{L}\leq L$,  we conclude that the PGD applied to $\wJ(\xb)$ has a better asymptotic convergence rate than standard GD applied to $J(\zb)$. 

To conclude this section, we present on the right of Fig. \eqref{fig:quadratic} a numerical study of the convergence of GD and PGD, with Armijo backtracking, for the strongly convex function
\begin{equation}\label{eq:strongly_convex}
J(\zb):=\text{Log}\left(\sum_{i=1}^n a_i e^{b_i x_i}\right)+ \frac{1}{2} \zb^\top D \zb,
\end{equation}
where $n=10^3$, $n_{el}=20$, $a_i=i$ $\forall i$, $b_i=10$, $1\leq i \leq n_{el}$, $b_i=1$, $i>n_{el}$, and $D$ is a diagonal matrix, with the first $n_{el}$ entries equal to $10^{-4}$ and the remaining ones equal to $10^{-2}$. The function is ill conditioned with respect to the first $n_{el}$ variables. By nonlinear eliminating them, the improvement of PGD over GD is impressive: GD requires 693 iterations and 3.74 seconds, while PGD needs 9 iterations and 0.25 seconds with a speed up of order 15. Notice that in this setting an Armijo backtracking line search is employed which, compared to the quadratic case, leads to further evaluations of the implicit map $h$ that is computed with an inner Newton method.

\section{Inexact variable reduction}
The gain due to the reduced number of PGD iterations may be null when the evaluation of $h(\cdot)$ is expensive. In such cases, one could consider an inexact variable reduction in which the set of nonlinear equations $\nabla_{\yb} J(\xb,\yb)=0$ are solved approximately with an iterative scheme. An example is the iterative procedure
\begin{equation}\label{eq:inexact_h}
\begin{split}
&\widehat{h}_N(\xb;\yb^{0}):=\yb^{N}, \\
&\yb^{\ell+1} = \yb^{\ell} - t_{\ell} D_{\ell}(\xb,\yb^{\ell})\nabla_{\yb} J(\xb,\yb^{\ell}), \text{ for $\ell=0,\dots,N-1$,}
\end{split}
\end{equation}
which performs $N$ steps of a (left) preconditioned GD (possibly Newton steps if $D_{\ell}(\xb,\yb^{\ell})=\left(\nabla^2_{yy} J(\xb,\yb^{\ell})\right)^{-1}$) on the optimization problem $\min_{\yb} J(\xb,\yb)$ starting from $\yb^{0}$ with $\xb$ fixed. Once a new iterate $\xb^{k+1}$ is computed, one can update the initial guess of \eqref{eq:inexact_h} by setting $\yb^{0}\leftarrow\widehat{h}_N(\xb^{k+1};\yb^{0})$.

To preserve consistency between the full and reduced optimization problems, the inexact elimination must satisfy $\widehat{h}_N(\xb^\star;\yb^\star)=\yb^\star$, that is, the iterative procedure, when starts from the optimal $\yb^\star$ associates to $\xb^\star$ again $\yb^\star$. This is trivially satisfied by \eqref{eq:inexact_h} since $\nabla_{\yb}J(\xb^\star,\yb^\star)=0$.
Incidentally, we emphasize that a special choice is $N=0$ leading to the constant map $\widehat{h}_0(\xb;\yb_0)=\yb_0$. We then have
\[\min_{\xb} \wJ(\xb)=\min_{\xb} J(\xb,\widehat{h}_0(\xb;\yb_0))=\min_{\xb} J(\xb,\yb_0),\]
and hence the minimization of $\wJ$ is equivalent to the minimization of $J$ with the variable $\yb$ freezed. To recover convergence to $(\xb^\star,\yb^\star)$, one may alternate the minimization of $J$ with respect to $\xb$ and $\yb$,
\begin{equation}\label{eq:alternating_minimization}
\xb^{k+1}=\text{argmin}_{\xb} J(\xb,\yb^k)\quad\text{and}\quad \yb^{k+1}=\text{argmin}_{\yb} J(\xb^{k+1};\yb),    
\end{equation}
recovering the popular alternating minimization methods (\cite{Bertsekas,Bezdek2002,Beck2013}), that  can thus be seen as particular instances of a reduced variable approach by choosing a specific inexact function $\widehat{h}(\cdot)$.

Notice further that an inexact evaluation of $h$ complicates the computation of the gradient of the $\wJ$ since 
\[\nabla \wJ(\xb^k)=\nabla_{\xb} J(\xb^k,\widehat{h}_N(\xb^k,\yb^{0}))+ \nabla_{\yb} J(\xb^k,\widehat{h}_N(\xb^k,\yb^{0}))\widehat{h}_N^\prime(\xb^k,\yb^{0}),\]
where the second term now does not cancel due to the inexactness of $\widehat{h}_N$, and $\widehat{h}_N^\prime$ may involve high order derivatives of $J$.
As a consequence, in our numerical implementation we perform inexact elimination using a variant of \eqref{eq:inexact_h} consisting of an inexact Newton method \cite{Dembo}. The iterative procedure is stopped when $\|\nabla_{\yb} J(\xb,\yb)\|<\text{Tol}$. The number of inner iteration depends then on the current tolerance, i.e., $N=N(\text{Tol})$. The initial tolerance is $10^{-3}$ and after each gradient step on $\wJ(\xb)$, $\text{Tol}$ is decreased by a factor $\rho:=0.5$. 
We then use the inexact gradient $\nabla \wJ(\xb^k)\approx \nabla_{\xb} J(\xb^k,\widehat{h}_{N}(\xb^k,\yb^{0,k}))$ as a descent direction, whose inexactness though is controlled by that of $\widehat{h}_N(\cdot)$ and vanishes in the limit for $k\rightarrow \infty$.
We briefly mention that we have also implemented the inexact elimination \eqref{eq:inexact_h} using $N$-fixed steps of gradient descent, and used the full gradient $\nabla \wJ$ as descent direction.
The results though were not satisfactory compared with the former procedure. 

Table \ref{Tab:Inexa} reports the number of iterations and computational times of GD and PGD with both exact and inexact elimination to minimize \eqref{eq:strongly_convex} for increasing values of $n_{el}$, that is the number of variables that are responsible for the ill-conditioning.
While PGD with exact elimination becomes inefficient in terms of computational times as $n_{el}$ grows, PGD with inexact elimination results very robust both in terms of iterations and computational time, and significantly outperforms GD in all test cases.
\begin{table}
\centering
\begin{tabular}{|c|c|c|c|c| } 
 \hline
 $n_{el}$ & 10 & 50 & 200 &400\\ 
  \hline
 GD & 713 (3.73) & 657 (3.73) & 558 (296) & 487 (2.52) \\  \hline
 PGD-Ex & 9 (0.30) & 9 (0.30) & 9 (0.67) & 10 (3.19) \\  \hline
 PGD-In & 9 (0.25) & 9 (0.29) & 9 (0.29) & 10 (0.30) \\ 
 \hline
\end{tabular}
\caption{Number of iterations and computational time in seconds for a standard GD, a right PGD with exact elimination and a right PGD with inexact elimination to minimize \eqref{eq:strongly_convex} up to a tolerance of $10^{-6}$ on the relative gradient norm.}\label{Tab:Inexa}
\end{table}
\section{Conclusions}
In this contribution, we argued that nonlinear elimination can be used to reduce the number of variables in optimization problems, resulting in reduced objective functions that can potentially be easier to minimize with classical optimization algorithms. Future efforts will focus on analyzing theoretically the convergence of inexact elimination procedures, and the connections with the popular alternating minimization methods. It is also relevant to develop general criteria to identify efficient decomposition of the optimization variables. Concerning computational aspects, we plan to investigate our approach in more realistic problems, such those appearing in risk-adverse PDE-constrained optimization \cite{CMG}.

\end{document}